\documentstyle[12pt]{article}
\def \Z{\hbox{$Z\hskip -5.2pt Z$}}

\def\sZ{\hbox{$\sc Z\hskip -4.2pt Z$}}
\def \Q{\hbox{$Q\hskip -6.3pt \vrule height 6pt depth 0pt\hskip 6pt$}}

\def \C{\hbox{$C\hskip -5pt \vrule height 6pt depth 0pt \hskip 6pt$}}

\def\qed{\hfill \hfill \ifhmode\unskip\nobreak\fi\ifmmode\ifinner
         \else\hskip5pt\fi\fi
 \hbox{\hskip5pt\vrule width4pt height6pt depth1.5pt\hskip 1 pt}}
\def\a{\alpha}
\def\b{\beta}
\def\d{\delta}
\def\D{\Delta}

\def\G{\Gamma}
\def\l{\lambda}

\def\si{\sigma}

\def\sc{\scriptstyle}
\def\ssc{\scriptscriptstyle}
\def\dis{\displaystyle}
\def\cl{\centerline}

\def\sF{\hbox{$\sc I\hskip -2.5pt F$}}
\def\nl{\newline}
\def\ol{\overline}

\def\nob#1{$\mbox{#1}$}

\def\rar{\rightarrow}

\def\F{\hbox{$I\hskip -4pt F$}}\def \Z{\hbox{$Z\hskip -5.2pt Z$}}
\def\AA{{\cal A}}
\def\DD{{\cal D}}

\def\Lra{\Leftrightarrow}
\def\bs{\backslash}
\def\hs{\hspace*}
\def\vs{\vspace*}

\def\rb{\raisebox}

\def\ra{\rangle}
\def\ptl{\partial}
\def\la{\langle}
\def\ni{\noindent}
\def\hi{\hangindent}
\def\ha{\hangafter}
\def\limit{\rb{-6pt}{\mbox{$^{\ \dis\rm lim}_{n\rar\infty\,}$}}\def\iso{\cong}}

\begin{document}
\def\TIT{Some representations of nongraded Lie algebras
of generalized Witt type}
\def\ABS
{In a paper by Xu, some simple Lie algebras of generalized Cartan
type were constructed, using the mixtures of grading operators and
down-grading operators. Among them, are the simple Lie algebras of
generalized Witt type, which are in general nongraded and have no
torus. In this paper, some representations of these simple Lie
algebras of generalized Witt type are presented.}
\def\cover{{\par\ \par\ \par\ni\hi5ex\ha1
{\bf Article title:}\nl \TIT
\par\ \par\ \par\ni\hi5ex\ha1
{\bf Correspondence address:}\nl
  Dr. Yucai Su \nl
  Department of Applied Mathematics \nl
  Shanghai Jiaotong University \nl
  1954 Huashan Road, Shanghai 200030 \nl
  China \nl
  Email: kfimmi@public1.sta.net.cn
\par\ \par\ \par\ni\hi5ex\ha1
{\bf Running title:}\nl
  higher rank Virasoro and super-Virasoro algebras
\par\ \par\ \par\ni\hi5ex\ha1
{\bf Abstract:}\nl \ABS
\par\ \par\ni\hi5ex\ha1
{\bf Keywords:}\nl\KEYW
\par
}\setcounter{page}{0}\pagebreak}
\cl{{\bf Some representations of nongraded Lie algebras of
generalized Witt type}\footnote{AMS Subject Classification -
Primary: 17B10, 17B65, 17B66, 17B68, 17B70 \nl\hs{3.5ex} Supported
by a grant from National Educational Department of China. }}
\par\vs{4pt}\cl{(appeared in {\it J.Alg.} {\bf246} (2001), 721-738)}
\par\ \vs{-7pt}\par
\cl{Yucai Su\footnote{Email: ycsu@sjtu.edu.cn}} \cl{\small\it
Department of Mathematics, Shanghai Jiaotong University, Shanghai
200030, China}
\par
\cl{and}
\par
\cl{Jianhua Zhou\footnote{Email: jhzhou@seu.edu.cn}} \cl{\small\it
Department of Applied Mathematics, Southeast University Nanjing
210096, China}
\par
{\small \ABS}
\par\
\vs{-7pt}
\par
\ni {\bf 1. INTRODUCTION}
\par
The four well-known series of infinite dimensional simple Lie \nob{algebras}
of Cartan type have played important roles in the structure theory
of Lie algebras. \nob{Generalizations} of the simple Lie algebras of Witt
type have been obtained by Kawamoto [Kaw],
Dokovic and Zhao [DZ1,DZ2], Xu [X1] and Zhao [Z2].
Passman [P] studied the Lie algebras $\AA\DD=\AA\otimes\DD$ of generalized
Witt type constructed from the pair of a commutative associative algebra
$\AA$ with an identity \nob{element} and its commutative derivation subalgebra
$\DD$ over a field $\F$ of arbitrary \nob{characteristic.}
Xu [X1] studied some of these simple Lie algebras of generalized Witt type
and other generalized Cartan types Lie algebras,
based on the pairs of the tensor algebra of the group
algebra of an additive subgroup of ${\F}^n$ with the polynomial algebra
in several variables and the subalgebra of commuting locally finite
derivations.  Su, Xu and Zhang [SXZ], Su and Xu [SX]
determined the structure spaces of the generalized simple Lie algebras of Witt
type and of special type constructed in [X1].
Su and Zhao [SZ2] determined the second cohomology group and
gave some representations of some Lie algebras of generalized Witt type.
\par
The Lie algebras associated with vertex algebras, the Lie algebras generated
by \nob{conformal} algebras are in general \nob{nongraded} and \nob{nonlinear}
Lie algebras [BXZ,X2]. The algebraic \nob{aspects}
of quantum field theory are the
representation theory of the Lie algebras generated by \nob{conformal}
algebras [K].
However, not much work has been known on the representation theory of the
\nob{nongraded} and nonlinear Lie algebras.
Su and Zhao [SZ3] defined the Lie \nob{algebras} of
Weyl type $\AA[\DD]=\AA\otimes\F[\DD]$, which are in general nongraded and
nonlinear, where $\AA$ is a \nob{commutative} associative
algebra with an identity element over a field $\F$ of arbitrary
\nob{characteristic}, and $\F[\DD]$ is the polynomial algebra  of
a commutative derivation subalgebra $\DD$ of $\AA$.
They also determined in [SZ4] the isomorphism classes and automorphism
groups of the associative and Lie algebras of Weyl type $\AA[\DD]$, where
$\F$ is a field of characteristic zero, $\DD$ is a commuting subalgebra of
locally finite but not locally nilpotent derivations of $\AA$.
\par
In this paper, we shall consider some representations of the
(in general, nongraded)
Lie algebras of generalized Witt type defined by Xu [X1] below.
\par
Denote by $\F$ an algebraically closed field of characteristic zero.
Let $\ell_1,\ell_2,\ell_3$ be three nonnegative integers such that
$\ell=\ell_1+\ell_2+\ell_3>0$.
For any $m,n\in\Z$, we use the notation
$$
\ol{m,n}=\{m,m+1,...,n\}\mbox{ if }m\le n,
\mbox{ \ \ or \ \ }\ol{m,n}=\emptyset\mbox{ if }m>n.
\eqno(1.1)$$
Take a nondegenerate additive subgroup $\G$ of $\F^{\ell_2+\ell_3}$ in the
sense that $\G$ contains an $\F$-basis of $\F^{\ell_2+\ell_3}$.
For convenience, elements $\a\in\G$ will be written as
$$\a=(\a_1,\a_2,...,\a_\ell)\mbox{ with }
\a_1=\a_2=...=\a_{\ell_1}=0.
\eqno(1.2)$$
For any $i\in\Z_+=\{0,1,2,...\}$, $p\in\ol{1,\ell}$, we denote
$$
i_{[p]}=(0,...,0,\rb{6pt}{\mbox{$^{^{\sc p}}_{\dis i}$}},0,...,0)\in\Z^\ell_+.
\eqno(1.3)$$
Elements in $\Z^\ell_+$ will be denoted by
$$
\vec i=(i_1,i_2,...,i_\ell),
\eqno(1.4)$$
while elements in $\Z^{\ell_1+\ell_2}_+$ will be denoted by (1.4) with
$i_{\ell_1+\ell_2+1}=...=i_{\ell-1}=i_\ell=0$.
\par
Let $\AA$ be the semi-group algebra
$\F[\G\times\Z_+^{\ell_1+\ell_2}]$
with \F-basis $\{x^\a t^{\vec i}\,|\,(\a,\vec i{\ssc\,})
\in\G\times\Z_+^{\ell_1+\ell_2}\}$ and the algebraic operation
$$
x^\a t^{\vec i}\cdot x^\b t^{\vec j}=
x^{\a+\b}t^{\vec i+\vec j}\,\,
\eqno(1.5)$$
for $(\a,\vec i{\ssc\,}),(\b,\vec j{\ssc\,})\in\G\times\Z_+^{\ell_1+\ell_2}.$
Denote the identity element $x^0t^0$ by 1.
Define the linear transformations (derivations)
$\{\ptl^-_1,\ptl^-_2,...,\ptl^-_{\ell_1+\ell_2},\,
\ptl^+_{\ell_1+1},\ptl^+_{\ell_1+2},...,\ptl^+_\ell\}$ on $\AA$ by
$$
\ptl^-_p(x^\a t^{\vec i})=i_p x^\a t^{\vec i-1_{[p]}},\ \
\ptl^+_q(x^\a t^{\vec i})=\a_q x^\a t^{\vec i},
\eqno(1.6)$$
for all $p\in\ol{1,\ell_1+\ell_2},\,q\in\ol{\ell_1+1,\ell}$, where
in the first equation, if $\vec i-1_{[p]}\notin\Z_+^{\ell_1+\ell_2}$ then
$i_p=0$ and the right-hand side is treated as zero.
The operators $\ptl_p^-$ are called {\it down-grading operators} and
$\ptl^+_q$ are {\it grading operators}.
Set
$$
\ptl_p=\ptl^-_p,\ \,\ptl_q=\ptl^-_q+\ptl^+_q,\ \,
\ptl_r=\ptl^+_r,
\eqno(1.7)$$
for $p\in\ol{1,\ell_1},\,q\in\ol{\ell_1+1,\ell_1+\ell_2},\,
r\in\ol{\ell_1+\ell_2+1,\ell}.$
Observe that $\ptl_p$ is locally nilpotent if $p\in\ol{1,\ell_1}$;
locally finite if $p\in\ol{\ell_1+1,\ell_1+\ell_2}$;
semi-simple if $p\in\ol{\ell_1+\ell_2+1,\ell}$.
Set $\DD=\sum_{p=1}^\ell\F\ptl_p$.
\par
It is proved in [SXZ] that the pairs $(\AA,\DD)$ constructed above
for different parameters
$(\ell_1,\ell_2,\ell_3;\G)$ enumerate all the pairs $(\AA,\DD)$
of a commutative associative algebra $\AA$ with an identity element and its
finite-dimensional commutative locally-finite derivation subalgebra $\DD$
such that $\AA$ is $\DD$-simple and
$\cap_{\ptl\in\DD}{\ssc\,}{\rm ker\ssc\,}\ptl=\F.$
\par
The \F-vector space
$$
W=W(\ell_1,\ell_2,\ell_3,\G)=\AA\otimes\DD=\sum_{p=1}^\ell\AA\ptl_p,
\eqno(1.8)$$
forms a Lie algebra under the usual bracket
$$
[x^\a t^{\vec i}\ptl_p,x^\b t^{\vec j}\ptl_q]=
x^{\a+\b}t^{\vec i+\vec j}(\b_p\ptl_q-\a_q\ptl_p)
+j_px^{\a+\b}t^{\vec i+\vec j-1_{[p]}}\ptl_q
-i_qx^{\a+\b}t^{\vec i+\vec j-1_{[q]}}\ptl_p,
\eqno(1.9)$$
for
$(\a,\vec i{\ssc\,}),(\b,\vec j{\ssc\,})\in\G\times\Z_+^{\ell_1+\ell_2},\,p,q\in\ol{1,\ell},$
which is called a {\it Lie Algebra of generalized Witt type} constructed
by Xu [X1].
In particular, $W(\ell,0,0,\{0\})$ is a classical Witt algebra
$W^+_{\ell}$ over the polynomial ring $\F[t_1,...,t_\ell]$, and
$W(0,0,\ell,\Z^{\ell})$ is a classical Witt algebra $W_\ell$
over the Laurent polynomial ring $\F[x_1^{\pm1},...,x_\ell^{\pm1}]$,
and more general,
$W(0,0,\ell,\G)$ is a generalized Witt algebra considered in
[Z1] and
$W(\ell_1,0,\ell_3,\G)$ is a generalized Witt algebra considered in
[DZ1,DZ2]. Also, observe that $W_1$ is the classical centerless Virasoro
algebra, $W(0,0,1,\G)$ is a centerless higher rank Virasoro algebra
(if $\G\subset\F$ is finitely generated) [PZ] or a centerless generalized
Virasoro algebra [SZ1] and $W(0,1,0,\G)$ is a nongraded centerless
Virasoro algebra [SL].
\par
Eswara \hfill Rao \hfill [E1,E2] \hfill constructed \hfill some \hfill representations \hfill of \hfill $W^+_\ell$, \hfill
Zhao \hfill [Z1] \hfill classified\nl
\nob{$W(0,0,\ell,\G)$-modules}
with weight multiplicity 1.
There is much more work on the Virasoro algebra and generalized Virasoro
algebras. Mathieu [M] classified Harish-Chandra modules over the Virasoro
algebra with a partial result also obtained in [S1].
Su generalized Mathieu's result to super-Virasoro algebras
and to higher rank Virasoro and super-Virasoro algebras [S2,S3,S4].
\par
Since $W(\ell_1,\ell_2,\ell_3,\G)$ in general does not have a torus,
there is no concept of weight modules. However, motivated by [SL],
we shall be able to consider the so-called {\it generalized weight modules}
in Section 2 and give a classification of generalized weight modules with
weight multiplicity 1 in Section 3. Our main result is Theorem 2.3.
\par
Finally, \hfill we \hfill shall \hfill remark \hfill that \hfill our \hfill classification \hfill also \hfill applies \hfill to \hfill an \hfill arbitrary \hfill
field \hfill of\nl
\nob{characteristic} zero as in Section 5 of [Z1].
\par\ \vs{-7pt}
\par\ni
{\bf 2. PRELIMINARY AND MAIN RESULT}
\par
The Lie algebra $W=W(\ell_1,\ell_2,\ell_3,\G)$ is in general nongraded in
the sense that it can not be graded so that all homogeneous spaces are
finite dimensional.
Since $W$ in general does not contain a toral subalgebra, there is no
concept of weight modules. However, we shall be able to consider the
modules defined below.
\par
A linear transformation $T$ of a $\F$-vector space $V$ is called
{\it locally finite} if
$$
{\rm dim}({\rm span}_{\sF}\{T^n(v)\,|\,n\in\Z_+\})<\infty,
\eqno(2.1)$$
for all $v\in V$. $T$ is called {\it semi-simple} if there exists
a \F-basis of ${\rm span}_{\sF}\{T^n(v)\,|\,n\in\Z_+\}$ consisting of
eigenvectors of $T$ for all $v\in V$; {\it locally nilpotent}
if for any $v\in V$, there exists $n\in\Z_+$ (may depend on $v$)
such that $T^n(v)=0$. A subspace $U$ of $End(V)$ is called {\it locally finite} or {\it locally nilpotent} if every element of $U$ is  {\it locally finite} or {\it locally nilpotent} on $V$ respectively.
\par
Denote
$$
\DD_1=\sum_{p=1}^{\ell_1}\F\ptl_p,\ \
\DD_2=\sum_{p=\ell_1+1}^{\ell_1+\ell_2}\F\ptl_p,\ \
\DD_3=\sum_{p=\ell_1+\ell_2+1}^\ell\F\ptl_p.
\eqno(2.2)$$
For convenience, an element $\b\in\F^{\ell_2+\ell_3}$ will be denoted as in
(1.2) by $\b=(\b_1,\b_2,...,\b_\ell)$
with $\b_1=\b_2=...=\b_{\ell_1}=0$.
Let $\ptl=\sum_{p=1}^\ell a_p\ptl_p\in\DD$, we define
$$
\la\ptl,\b\ra=\b(\ptl)=\sum_{p=\ell_1+1}^\ell a_p\b_p.
\eqno(2.3)
$$
By (1.6), (1.7), we see that $\rm ad\,\DD_1$ is locally
nilpotent and $\rm ad\,\DD$ is locally finite on $W$.
Thus, we can consider $W$-modules $V$ such that $\DD_1$ and $\DD$
act respectively locally nilpotently and locally finitely on $V$.
For such a module $V$, we have
$$
V=\bigoplus_{\b\in\sF^{\ell_2+\ell_3}}V_\b,\ \
V_\b=\{v\in V\,|\,(\ptl-\b(\ptl))^n(v)=0\mbox{ for $\ptl\in\DD$ and
some }n\in\Z_+\}.
\eqno(2.4)$$
Let $\b\in\F^{\ell_2+\ell_3}$. We set
$$
V_\b^{(n)}=\{v\in V_\b\,|\,(\ptl-\b(\ptl))^{n+1}(v)=0\mbox{ for }\ptl\in\DD\},
\eqno(2.5)$$
for $n\ge0$; and set
$$
\ol V_\b^{(0)}=V_\b^{(0)},\ \
\ol V_\b^{(n)}=V_\b^{(n)}/V_\b^{(n-1)},
\eqno(2.6)$$
for $n\ge1$. Clearly,
$$
V_\b=\bigcup_{n\in\sZ_+}V_\b^{(n)}\mbox{ \ and \ }
V_\b\ne\{0\}\ \Lra\ V_\b^{(0)}\ne\{0\}.
\eqno(2.7)$$
For $n\ge0$, we set
$$
V^{(n)}=\bigoplus_{\b\in\sF^{\ell_2+\ell_3}}V_\b^{(n)},\ \
\ol V^{(n)}=\bigoplus_{\b\in\sF^{\ell_2+\ell_3}}\ol V_\b^{(n)}.
\eqno(2.8)$$
The above discussion leads us to the following definition.
\par
{\bf Definition 2.1}. {\it A module $V$ over $W$ is called a generalized
weight module if $V$ has a \nob{decomposition} (2.4).
For $\b\in\F^{\ell_2+\ell_3}$, if
$V_\b\ne\{0\}$, then $V_\b$ is called the generalized weight space with
weight $\b$ and $V_\b^{(0)}$ is called the weight space with weight $\b$.
For an \nob{indecomposable} \nob{generalized} weight module $V$ over $W$,
it is called a module of the intermediate series if
${\rm dim\ssc\,}V_\b^{(0)}\le1$ for all $\b\in\F^{\ell_2+\ell_3}$;
it is called a uniformly bounded module
if there exists a nonnegative integer $N$ such that
${\rm dim\ssc\,}V_\b^{(0)}\le N$ for all $\b\in\F^{\ell_2+\ell_3}$.}
\qed\par
The aim of this paper is to classify
all modules of the intermediate series
over $W$. We shall assume that $\ell_1+\ell_2\ge1$ since the case for
$\ell_1+\ell_2=0$ has been considered by Zhao [Z1].
We shall also assume that $\ell_2+\ell_3\ge1$ since
$W$ is a classical Witt algebra if $\ell_2+\ell_3=0$.
Note that $W$ is a
centerless nongraded Virasoro algebra considered in [SL] if $\ell=\ell_2=1$.
For convenience to
the reader, we state below the main classification theorem in [Z1].
\par
{\bf Theorem 2.2}. {\it (1) A module of the intermediate series over
$W(0,0,\ell,\G)$ with $\ell\ge2$ is a quotient module of $A_{\a,b},A(\b)$
or $B(\b)$ for
some suitable $\a,\b\in\F^\ell,b\in\F$, where $A_{\a,b},A(\b),B(\b)$ all
have $\F$-basis $\{v_\mu\,|\,\mu\in\G\}$ with the following actions:
$$
\matrix{
A_{\a,b}:\hfill\!\!\!\!&
(x^\mu\ptl)v_\nu=\la\ptl,\a+\nu+b\mu\ra v_{\mu+\nu},
\vs{4pt}\hfill\cr
A(\b){\ssc\!}:\hfill\!\!\!\!&
(x^\mu\ptl)v_\nu=\la\ptl,\nu\ra v_{\mu+\nu},
\ \nu,\mu+\nu\ne0,\ (x^\mu\ptl)v_{-\mu}=0,\
(x^\mu\ptl)v_0=\la\ptl,\mu+\b\ra v_{\mu},
\vs{4pt}\hfill\cr
B(\b){\ssc\!}:\hfill\!\!\!\!&
(x^\mu\ptl)v_\nu\!=\!\la\ptl,\mu\!+\!\nu\ra v_{\mu+\nu},
\ \nu,\mu\!+\!\nu\!\ne\!0,\, (x^\mu\ptl)v_0\!=\!0,\,
(x^\mu\ptl)v_{-\mu}\!=\!-\la\ptl,\mu\!+\!\b\ra v_0,
\hfill\cr
}
\eqno\matrix{\hfill(2.9)\vs{4pt}\cr(2.10)\vs{4pt}\cr(2.11)\cr}$$
for $\mu,\nu\in\G,\ptl\in\DD$.
\vs{-3pt}\par
(2) $A_{\a,b}$ is simple if and only if $\a\notin\G$ or $b\ne0,1$.
\vs{-3pt}\par
(3) each of $A_{0,0},A_{0,1},A(\b),B(\b)$ has two composition factors:
one is 1-dimensional (trivial)
and the other is  simple. The simple ones are denoted  respectively by
$A'_{0,0},A'_{0,1},A'(\b),B'(\b)$. Then
$A'_{0,1},A'(\b)$ are simple submodules spanned by all basis element except
$v_0$.
$A'_{0,0},B'(\b)$ are simple quotient modules modulo the trivial submodule
spanned by $v_0$. We have $A'_{0,0}\cong A'(\b),A'_{0,1}\cong B'(\b)$.
}
\qed\par
In this paper, we shall use the convention that
if a notation which is not defined but techniquely appears in an expression,
we always
treat it as zero; for example, in the following theorem, we treat
$v_{\mu,\vec i}$ as zero if $\vec i\notin \Z_+^{\ell_1+\ell_2}$.
\par
The main result of this paper is the following theorem.
\par
{\bf Theorem 2.3}. {\it (1) Suppose $\ell_1+\ell_2\ge1, \ell_2+\ell_3\ge1$.
A module $V$ of the intermediate series over
$W(\ell_1,\ell_2,\ell_3,\G)$ is
a subquotient module of $A_{\a,b}$ for some suitable
$\a\in\F^{\ell_2+\ell_3},b\in\F$, where $A_{\a,b}$ is a module
with \F-basis $\{v_{\mu,\vec i}\,|\,\mu\in\G,\vec i\in\Z_+^{\ell_1+\ell_2}\}$
and the following action
$$
A_{\a,b}:
(x^\mu t^{\vec i}\ptl_p)v_{\nu,\vec j}
=(\a_p+\nu_p+b\mu_p)v_{\mu+\nu,\vec i+\vec j}
+(j_p+bi_p)v_{\mu+\nu,\vec i+\vec j-1_{[p]}},
\eqno(2.12)$$
for $\mu,\nu\in\G,\,\vec i,\vec j\in\Z_+^{\ell_1+\ell_2},\,p\in\ol{1,\ell}$,
where if $i_p+j_p-1<0$ then $i_p=j_p=0$ and so $j_p+bi_p=0$
thus the second term of the right-hand side does not occur.
\vs{-3pt}\par
(2) $A_{\a,b}$ is simple if and only if $\a\notin\G$ or $b\ne0$.
\vs{-3pt}\par
(3) $A_{0,0}$ has two composition factors:
the trivial submodule $\F v_{0,0}$ and the simple quotient module
$A'_{0,0}=A_{0,0}/\F v_{0,0}$. However, if $\ell_1+\ell_2\ge2$, then
$V'=A'_{0,0}$ is not a module of the intermediate series since $V'{}_0^{(0)}
={\rm span}_{\sF}\{v_{0,1_{[p]}}\,|\,p\in\ol{1,\ell_1+\ell_2}\}$ has dimension
$\ell_1+\ell_2>1$.
\vs{-3pt}\par
(4) The possible isomorphisms between modules of the intermediate series are
the following:
(i) $A_{\a,b}\cong A_{\b,d}\Lra \a-\b\in\G,b=d$.
(ii) $A'_{0,0}\cong A_{0,1}$ if $\ell=\ell_2=1$.
}
\par
We would like to remark that unlike the case $\ell_1+\ell_2=0$, we do not
have modules of types $A(\a),B(\a)$ when $\ell_1+\ell_2>0$.
It is also interesting to see that $A_{0,1}$ is simple, which is not the
case for the Witt algebras (cf. Theorem 2.2(2)) and that
$A'_{0,0}\cong A_{0,1}$ if $\ell=\ell_2=1$.
\par\ \vs{-7pt}\par\ni
{\bf 3. PROOF OF THEOREM 2.3}
\par
Suppose that $V$ is a simple module of the intermediate series over
$W=W(\ell_1,\ell_2,\ell_3,\G)$, where $\ell_1+\ell_2\ge1,\ell_2+\ell_3\ge1$.
We shall further suppose $\ell\ge2$ since $\ell=\ell_2=1$ has been considered
in [SL].
For any $\a\in\F^{\ell_2+\ell_3}$,
let $V(\a)=\oplus_{\mu\in\G}V_{\a+\mu}$.
Then obviously, $V(\a)$ is a submodule of $V$ and $V$ is a direct sum of
different $V(\a)$. Since $V$ is simple, we must have $V=V(\a)$ for
some $\a\in\F^{\ell_2+\ell_3}$.
\par
Define a total order on $\Z_+^{\ell_1+\ell_2}$ as follows:
$$
\vec i<\vec j\Lra|\vec i|<|\vec j|
\mbox{ or }|\vec i|=|\vec j|
\mbox{ but }\exists\,p\mbox{ such that }i_p<j_p\mbox{ and }
i_q=j_q\mbox{ for }q\in\ol{1,p-1},
\eqno(3.1)$$
for $\vec i,\vec j\in\Z_+^{\ell_1+\ell_2}$,
where $|\vec i|=\sum_{p=1}^{\ell_1+\ell_2}i_p$ is the {\it level}
of $\vec i$.
\par
We shall prove Theorem 2.3(1) in two cases.
\par
{\it Case 1}. $\a\notin\G$.
\par
Observe that by (1.9), ${\rm span}_{\sF}\{x^\mu\ptl\,|\,\mu\in\G,\ptl\in
\DD_2+\DD_3\}$ is the Witt algebra $W(0,0,\ell_2+\ell_3,\G)$.
Note that $V^{(0)}=\oplus_{\mu\in \G}V_{\a+\mu}^{(0)}$ is a
$W(0,0,\ell_2+\ell_3,\G)$-submodule of $V$. By Theorem 2.2, it is a simple
\nob{$W(0,0,\ell_2+\ell_3,\G)$-module} of the intermediate series of type $A_{\a,b}$ for some
$b\in\F$. Thus there
exists a \F-basis $\{v_{\mu,0}\,|\,\mu\in\G\}$ of $V^{(0)}$ such that
$$
(x^\mu\ptl)v_{\nu,0}=\la\ptl,\a+\nu+b\mu\ra v_{\mu+\nu,0},
\eqno(3.2)$$
for $\mu,\nu\in\G,\ptl\in\DD_2+\DD_3$.
We shall show that (3.2) also holds for $\ptl\in\DD_1$
(It is possible that $\ell_1=0$, in this case, the arguments from here to
(3.7) are trivial). Thus let $\ptl
=\ptl_p$ with $p\in\ol{1,\ell_1}$. Then
$$(\ptl_q-(\a_q+\nu_q))
(t^{1_{[p]}}\ptl_p)v_{\nu,0}=
\d_{p,q}\ptl_pv_{\nu,0}
+(t^{1_{[p]}}\ptl_p)(\ptl_q-(\a_q+\nu_q))v_{\nu,0}=0,
\eqno(3.3)$$
for $q\in\ol{\ell_1+1,\ell}$. Thus $(t^{1_{[p]}}\ptl_p)v_{\nu,0}\in V_{\a+\nu}^{(0)}$
and so
$$(t^{1_{[p]}}\ptl_p)v_{\nu,0}=c_{p,\nu}v_{\nu,0}
\mbox{ for some }c_{p,\nu}\in\F.
\eqno(3.4)$$
Similarly, we have
$$
(x^\mu\ptl_p)v_{\nu,0}=d_{p,\mu,\nu} v_{\mu+\nu,0}
\mbox{ for some }d_{p,\mu,\nu}\in\F.
\eqno(3.5)$$
Applying ${\rm ad\ssc\,}x^\mu\ptl_q$ for $q\in\ol{\ell_1+1,\ell}$ to (3.4), using
(3.2), we obtain that
$c_{p,\nu}=c_{p,\mu+\nu}$ if $\a+\nu+b\mu\ne0$, from this it is easily
derived that $c_{p,\nu}=c_p$ does not depend on $\nu\in\G$.
Now applying ${\rm ad\ssc\,}x^\mu\ptl_p$ to (3.4), we obtain
$$
(t^{1_{[p]}}\ptl_p)(x^\mu\ptl_p)v_{\nu,0}=(c_p-1)(x^\mu\ptl_p)v_{\nu,0}.
\eqno(3.6)$$
However, using (3.4), (3.5), $(x^\mu\ptl_p)v_{\nu,0}$, if not zero,
is an eigenvector of $t^{1_{[p]}}\ptl_p$ corresponding to eigenvalue
$c_{p,\mu+\nu}=c_p$. Comparing this with (3.6) gives
$$
(x^\mu\ptl_p)v_{\nu,0}=0,
\eqno(3.7)$$
for $p\in\ol{1,\ell_1}$, which coincides with (3.2).
\par
For any
$\mu\in\G,0\ne\vec i\in\Z_+^{\ell_1+\ell_2}$, we define $v_{\mu,\vec i}$
inductively as follows. Take
$$
p'=p'_{\vec i}={\rm min}\{p\in\ol{1,\ell_1+\ell_2}\,|\,i_p\ne0\},\ \
p''=p''_\mu={\rm min}\{p\in\ol{\ell_1+1,\ell}\,|\,\a_p+\mu_p\ne0\}.
\eqno(3.8)$$
Since $\a\notin\G$, for any $\mu\in\G$, we have $\a+\mu\ne0$, i.e.,
there exists $p\in\ol{\ell_1+1,\ell}$ with $\a_p+\mu_p\ne0$. Thus $p''$
in (3.8) is uniquely defined.
Define
$$v_{\mu,\vec i}=
(\a_{p''}+\mu_{p''})^{-1}
((t^{1_{[p']}}\ptl_{p''})v_{\mu,\vec i-1_{[p']}}-
(i_{p''}-\d_{p',p''}+\d_{p',p''}b)v_{\mu,\vec i-1_{[p'']}}).
\eqno(3.9)$$
\par
{\bf Claim 1}. For $p\in\ol{1,\ell},(\mu,\vec i{\ssc\,})\in\G\times
\Z_+^{\ell_1+\ell_2}$, we have
$$
\ptl_p v_{\mu,\vec i}
=(\a_p+\mu_p)v_{\mu,\vec i}+i_p v_{\mu,\vec i-1_{[p]}}.
\eqno(3.10)$$
\par
If $\vec i=0$, (3.10) follows from (3.2) and (3.7). Suppose $\vec i>0$.
By (3.9) and the inductive assumption on the order of $\vec i$, (3.10) is equivalent to
$$
\matrix{
(\a_{p''}+\mu_{p''})^{-1}
(\d_{p,p'}\ptl_{p''}v_{\mu,\vec i-1_{[p']}}
+(i_p-\d_{p,p'})(t^{1_{[p']}}\ptl_{p''})v_{\mu,\vec i-1_{[p]}-1_{[p']}}
\vs{4pt}\hfill\cr
-(i_{p''}-\d_{p',p''}+\d_{p',p''}b)
(i_p-\d_{p,p''})v_{\mu,\vec i-1_{[p]}-1_{[p'']}})
\vs{4pt}\hfill\cr
=i_pv_{\mu,\vec i-1_{[p]}}.
\hfill\cr}
\eqno(3.11)$$
We shall verify (3.11) case by case as follows. {\it Case (i)}:
$p\ne p',p''$. In this case, since $p'_{\vec i-1_{[p]}}=p'$,
(3.11) is precisely the definition (3.9) with $\vec i$ replaced
by $\vec i-1_{[p]}$. {\it Case (ii)}: $p=p'\ne p''$. Then (3.11) becomes
$$
\matrix{
(\a_{p''}+\mu_{p''})^{-1}
((\a_{p''}+\mu_{p''})v_{\mu,\vec i-1_{[p']}}
+i_{p''}v_{\mu,\vec i-1_{[p']}-1_{[p'']}}
+(i_{p'}-1)(t^{1_{[p']}}\ptl_{p''})v_{\mu,\vec i-2_{[p']}}
\vs{4pt}\hfill\cr
-i_{p''}i_{p'}v_{\mu,\vec i-1_{[p']}-1_{[p'']}})
\vs{4pt}\hfill\cr
=i_{p'}v_{\mu,\vec i-1_{[p']}}.
\hfill\cr}
\eqno(3.12)$$
Since $i_{p'}>0$ by (3.8), if $i_{p'}=1$, then (3.12) holds trivially;
if $i_{p'}>1$, then $p'_{\vec i-1_{[p']}}=p'$ and (3.12) holds by definition
(3.9) with $\vec i$ replaced by $\vec i-1_{[p']}$.
The proofs for {\it Case (iii)}: $p=p''\ne p'$ and
{\it Case (iv)}: $p=p'=p''$ are exactly analogous.
\par
{\bf Claim 2}. For $p\in\ol{1,\ell_1+\ell_2},\,q\in\ol{1,\ell},\,
\mu,\l\in\G$, we have
$$
\matrix{
(x^\mu\ptl_q)v_{\l,1_{[p]}}\hfill&\!\!\!\!=&\!\!\!\!
(\a_q+\l_q+b\mu_q)v_{\mu+\l,1_{[p]}}+\d_{p,q}v_{\mu+\l,0},
\vs{4pt}\hfill\cr
(x^\mu t^{1_{[p]}}\ptl_q)v_{\l,0}\hfill&\!\!\!\!=&\!\!\!\!
(\a_q+\l_q+b\mu_q)v_{\mu+\l,1_{[p]}}+\d_{p,q}bv_{\mu+\l,0}.
\hfill\cr}
\eqno\matrix{(3.13)\vs{4pt}\cr(3.14)\cr}$$
\par
Let $p''=p''_\l$ be as defined in (3.8). Then $p''\ge\ell_1+1$. If
$p''>\ell_1+1$, using inductive assumption, we can suppose that
(3.13) holds for all
$\l\in\G$ with $p''_\l<p''$.
\par
As in the proof of (3.4), it is straightforward to verify
$$
(x^\mu\ptl_q)v_{\l,1_{[p]}}-(\a_q+\l_q+b\mu_q)v_{\mu+\l,1_{[p]}},
(x^\mu t^{1_{[p]}}\ptl_q)v_{\l,0}-(\a_q+\l_q+b\mu_q)v_{\mu+\l,1_{[p]}}
\in V_{\a+\mu+\l}^{(0)}.
\eqno(3.15)$$
Thus we can suppose
$$
\matrix{
(x^\mu\ptl_q)v_{\l,1_{[p]}}\hfill&\!\!\!\!=&\!\!\!\!
(\a_q+\l_q+b\mu_q)v_{\mu+\l,1_{[p]}}+
(\d_{p,q}+d_{p,q}^{\mu,\l})v_{\mu+\l,0},
\vs{4pt}\hfill\cr
(x^\mu t^{1_{[p]}}\ptl_q)v_{\l,0}\hfill&\!\!\!\!=&\!\!\!\!
(\a_q+\l_q+b\mu_q)v_{\mu+\l,1_{[p]}}
+(\d_{p,q}b+e_{p,q}^{\mu,\l})v_{\mu+\l,0},
\hfill\cr}
\eqno\matrix{(3.16)\vs{4pt}\cr(3.17)\cr}$$
for some $d_{p,q}^{\mu,\l},e_{p,q}^{\mu,\l}\in\F$. By (3.9), (3.10), we have
$$
d_{p,q}^{0,\l}=0,\ e_{p,p''}^{0,\l}=0.
\eqno(3.18)$$
We shall be careful that the second equation of (3.18) only holds for $p''$
not for general $q$ by definition (3.9) and (3.17).
Using
(3.16), (3.17) in
$$
\matrix{
[x^\mu\ptl_q,x^\nu\ptl_r]v_{\l,1_{[p]}}=
x^{\mu+\nu}(\nu_q\ptl_r-\mu_r\ptl_q)v_{\l,1_{[p]}},
\vs{4pt}\hfill\cr
[x^\mu\ptl_q,x^\nu t^{1_{[p]}}\ptl_r]v_{\l,0}=
(x^{\mu+\nu}t^{1_{[p]}}(\nu_q\ptl_r-\mu_r\ptl_q)+\d_{p,q}x^{\mu+\nu}\ptl_r)
v_{\l,0},
\hfill\cr}
\eqno(3.19)$$
we obtain (in the following, for convenience, for any $\mu\in\G$, we denote
$\ol\mu=\a+\mu$)
$$
\matrix{
(\ol\l_r+b\nu_r)d_{p,q}^{\mu,\nu+\l}
+(\ol\l_q+\nu_q+b\mu_q)d_{p,r}^{\nu,\l}
-(\ol\l_q+b\mu_q)d_{p,r}^{\nu,\mu+\l}
-(\ol\l_r+\mu_r+b\nu_r)d_{p,q}^{\mu,\l}
\vs{4pt}\hfill\cr
-\nu_q d_{p,r}^{\mu+\nu,\l}+\mu_rd_{p,q}^{\mu+\nu,\l}=0,
\vs{8pt}\hfill\cr
(\ol\l_r+b\nu_r)d_{p,q}^{\mu,\nu+\l}
+(\ol\l_q+\nu_q+b\mu_q)e_{p,r}^{\nu,\l}-(\ol\l_q+b\mu_q)e_{p,r}^{\nu,\mu+\l}
-\nu_qe_{p,r}^{\mu+\nu,\l}+\mu_re_{p,q}^{\mu+\nu,\l}=0.
\hfill\cr}
\eqno\matrix{(3.20)\vs{10pt}\cr\rb{-10pt}{(3.21)}\cr}$$
\par
First assume that $q\in\ol{1,\ell_1}$. Then the second, third and
fifth terms of (3.20) vanish, so do the second, third, fourth terms of (3.21).
Thus
$$
\matrix{
(\ol\l_r+b\nu_r)d_{p,q}^{\mu,\nu+\l}
-(\ol\l_r+\mu_r+b\nu_r)d_{p,q}^{\mu,\l}
+\mu_rd_{p,q}^{\mu+\nu,\l}=0,
\vs{4pt}\hfill\cr
(\ol\l_r+b\nu_r)d_{p,q}^{\mu,\nu+\l}+\mu_re_{p,q}^{\mu+\nu,\l}=0,
\hfill\cr}
\eqno\matrix{(3.22)\vs{4pt}\cr(3.23)\cr}$$
for $q\in\ol{1,\ell_1}$.
Let $\nu=0$ in (3.23), we obtain
$$
\ol\l_r d_{p,q}^{\mu,\l}+\mu_r e_{p,q}^{\mu,\l}=0,
\eqno(3.24)$$
for $q\in\ol{1,\ell_1}$.
Subtract (3.23) from (3.22) and making use of
(3.24), we obtain
$$
(\mu_r+\nu_r)(\ol\l_r+\mu_r+b\nu_r)d_{p,q}^{\mu,\l}
-\mu_r(\ol\l_r+\mu_r+\nu_r)d_{p,q}^{\mu+\nu,\l}
=0,
\eqno(3.25)$$
for $q\in\ol{1,\ell_1}$.
Replacing $\nu$ by $\mu$, and replacing $\mu$ by $2\mu$ and $\nu$ by $-\mu$
in (3.25), we obtain respectively
$$
\matrix{
2\mu_r(\ol\l_r+(1+b)\mu_r)d_{p,q}^{\mu,\l}
-\mu_r(\ol\l_r+2\mu_r)d_{p,q}^{2\mu,\l}
=0,
\vs{4pt}\hfill\cr
\mu_r(\ol\l_r+(2-b)\mu_r)d_{p,q}^{2\mu,\l}
-2\mu_r(\ol\l_r+\mu_r)d_{p,q}^{\mu,\l}
=0.
\hfill\cr}
\eqno(3.26)$$
By (3.18), we can suppose $\mu\ne0$. Choose $r\in\ol{\ell_1+1,\ell}$ with
$\mu_r\ne0$. If $b\ne0,1$, then the calculation of the determinant of the coefficients of
(3.26) shows that $d_{p,q}^{\mu,\l}=0$ and so $e_{p,q}^{\mu,\l}=0$ by (3.23)
for all $p\in\ol{1,\ell_1+\ell_2},q\in\ol{1,\ell_1},\mu,\l\in\G$.
\par
Assume that $b=0$ (the proof for $b=1$ is exactly analogous).
Replacing $\nu$ by $\nu-\mu$ in (3.25) shows that
$$
d_{p,q}^{\mu,\l}={\mu_r\over\ol\l_r+\mu_r}c_{p,q},\
e_{p,q}^{\mu,\l}={-\ol\l_r\over\ol\l_r+\mu_r}c_{p,q},\mbox{ for some }
c_{p,q}\in\F,
\eqno(3.27)$$
if $\l_r,\mu_r,\ol\l_r+\mu_r\ne0$, where the second equation follows from
(3.24). Using this and (3.23), we obtain that (3.27) holds for all $\mu,\l$
with $\ol\l_r+\mu_r\ne0$. In particular, $e_{p,q}^{0,\l}=c_{p,q}$.
As in the proof of (3.7), we have $(t^{1_{[p]}}\ptl_q)v_{\l,0}=0$
if $p\ne q$ and $q\in\ol{1,\ell_1}$, i.e., $c_{p,q}=0$ if
$p\ne q$ and $q\in\ol{1,\ell_1}$. Suppose $p=q\in\ol{1,\ell_1}$. Then by (3.17),
(3.27) and using $[t^{1_{[q]}}\ptl_q,t^{i_{[q]}}\ptl_q]=(i-1)t^{i_{[q]}}\ptl_q$,
we can prove by induction on $i$ that
$(t^{i_{[q]}}\ptl_q)v_{\l,0}=-ic_{q,q}v_{\l,i_{[q]}-1_{[q]}}$. Applying
${\rm ad\ssc\,}x^\mu\ptl_q$ to this gives that
$(x^\mu t^{i_{[q]}-1_{[q]}}\ptl_q)v_{\l,0}=
-c_{q,q}(x^\mu\ptl_q)v_{\l,i_{[q]}-1_{[q]}}$.
In particular, taking $i=2$ and comparing this with (3.17), (3.27), we obtain
that $c_{q,q}=0$. This proves that (3.13), (3.14) hold for
$q\in\ol{1,\ell_1}$.
\par
Now assume that $q\in\ol{\ell_1+1,\ell}$.
Our strategy is to prove that $d^{\mu,\nu}_{p,q}=0$
under some conditions on $\mu,\nu$ and then to obtain
the result in general. Thus, first we assume that $\l$ satisfies
$$
\ol\l_q\ne0,\l_q\ne0\mbox{ for all }q\in\ol{\ell_1+1,\ell}.
\eqno(3.28)$$
\par
Suppose that $\mu\in\G$ satisfies
$$
\mu_q\ne0,\,\forall\,q\in\ol{\ell_1+1,\ell},\mbox{ and }
\l_{\ell_1+1}\pm\mu_{\ell_1+1}\ne0\mbox{ and so }
p''_{\l\pm\mu}=\ell_1+1\le p''.
\eqno(3.29)$$
Then by taking $\nu=0$ and $r=p''$ in (3.21), using (3.18)
or the inductive assumption on $p''$, we obtain that
$$
\ol\l_{p''}d_{p,q}^{\mu,\l}+\mu_{p''}e_{p,q}^{\mu,\l}=0,
\eqno(3.30)$$
holds under condition (3.29).
Let $\nu=\mu,r=q$ in (3.21), using (3.30), we obtain that
$$
(\ol\l_q+b\mu_q)(\ol\l_{p''}+2\mu_{p''})d_{p,q}^{\mu,\mu+\l}
-\ol\l_{p''}(\ol\l_q+(b+1)\mu_q)d_{p,q}^{\mu,\l}=0,
\eqno(3.31)$$
holds under condition (3.29).
Let $\nu=-\mu,r=q$ in (3.20), we obtain
$$
(\ol\l_q-b\mu_q)d_{p,q}^{\mu,\l-\mu}
+(\ol\l_q+(b-1)\mu_q)d_{p,q}^{-\mu,\l}
-(\ol\l_q+b\mu_q)d_{p,q}^{-\mu,\mu+\l}
-(\ol\l_q+(1-b)\mu_q)d_{p,q}^{\mu,\l}=0.
\eqno(3.32)$$
Replacing $\mu$ by $-\mu$ and $\l$ by $\l+\mu$ in (3.31), we obtain
$$
(\ol\l_q+(1-b)\mu_q)(\ol\l_{p''}-\mu_{p''})d^{-\mu,\l}_{p,q}
-(\ol\l_{p''}+\mu_{p''})(\ol\l_q-b\mu_q)d^{-\mu,\mu+\l}_{p,q}=0.
\eqno(3.33)$$
Replacing $\l$ by $\l-\mu$ in (3.31), we obtain
$$
(\ol\l_q+(b-1)\mu_q)(\ol\l_{p''}+\mu_{p''})d^{\mu,\l}_{p,q}
-(\ol\l_{p''}-\mu_{p''})(\ol\l_q+b\mu_q)d^{\mu,\l-\mu}_{p,q}=0.
\eqno(3.34)$$
Let $\tau=(\ol\l_{p''}^2-\mu_{p''}^2)(\ol\l_q^2-(b\mu_q)^2)$.
Multiplying (3.32) by $\tau$ and using (3.33), (3.34), we have
$$
\matrix{
(\ol\l_{p''}+\mu_{p''})^2(\ol\l_q-b\mu_q)^2(\ol\l_q+(b-1)\mu_q)
d^{\mu,\l}_{p,q}+(\ol\l_q+(b-1)\mu_q)\tau d^{-\mu,\l}_{p,q}
\vs{4pt}\hfill\cr
-(\ol\l_q+(1-b)\mu_q)\tau d^{\mu,\l}_{p,q}-(\ol\l_q+b\mu_q)^2(\ol\l_{p''}
-\mu_{p''})^2(\ol\l_q+(1-b)\mu_q)d^{-\mu,\l}_{p,q}=0.
\hfill\cr}
\eqno(3.35)$$
i.e.
$$
(\ol\l_{p''}+\mu_{p''})(\ol\l_q-b\mu_q)\si d^{\mu,\l}_{p,q}+(\ol\l_{p''}-\mu
_{p''})(\ol\l_q+b\mu_q)\si d^{-\mu,\l}_{p,q}=0,
\eqno(3.36)$$
where
$$
\matrix{
\si(\mu)\!\!\!\!\!&
=\!(\ol\l_{p''}\!+\!\mu_{p''})(\ol\l_q\!-\!b\mu_q)
(\ol\l_q+(b\!-\!1)\mu_q)-(\ol\l_q+(1\!-\!b)\mu_q)
(\ol\l_{p''}\!-\!\mu_{p''})(\ol\l_q\!+\!b\mu_q)
\vs{4pt}\hfill\cr
&=2((b-2)\ol\l_{p''}\ol\l_q\mu_q+\ol\l_q^2\mu_{p''}+b(1-b)\mu_{p''}\mu_q^2).
\hfill\cr}
\eqno(3.37)$$
Thus if we assume that $\si(\mu)\ne0$, then we obtain
$$
(\ol\l_{p''}+\mu_{p''})(\ol\l_q-b\mu_q) d^{\mu,\l}_{p,q}+(\ol\l_{p''}-\mu_{p
''})(\ol\l_q+b\mu_q) d^{-\mu,\l}_{p,q}=0.
\eqno(3.38)$$
\par
Assume that
$$
(\ol\l_q+(b+n)\mu_q)(\ol\l_{p''}+(n+2)\mu_{p''})\ne0
\mbox{ for all }n\in\Z_+=\{0,1,2,...\}.
\eqno(3.39)$$
Then using (3.31) recursively, we obtain
$$
d_{p,q}^{\mu,n\mu+\l}=
{\ol\l_{p''}(\ol\l_{p''}+\mu_{p''})(\ol\l_q+(b+n)\mu_q)\over
(\ol\l_q+b\mu_q)(\ol\l_{p''}+n\mu_{p''})(\ol\l_{p''}+(n+1)\mu_{p''})}
d_{p,q}^{\mu,\l}\mbox{ for $n\in\Z_+$.}
\eqno(3.40)$$
Thus
$$
\limit
nd_{p,q}^{\mu,n\mu+\l}=
{\ol\l_{p''}(\ol\l_{p''}+\mu_{p''})\mu_q\over
\mu^2_{p''}(\ol\l_q+b\mu_q)}
d_{p,q}^{\mu,\l}\mbox{ \ for \ $n\in\Z_+$.}
\eqno(3.41)$$
(We shall remark that since we are encountering only rational functions
on $n$, and when $\mu,\l$, etc.~are fixed, we can regard functions as
defined in some extension field of the rational field $\Q$.
Thus the limit has meaning.)
In (3.20), setting $r=q$ and setting
$\mu,\nu,\l$ to be $2\mu,-\mu,2n\mu+\l\in\G$ respectively, we obtain
$$
\matrix{
 (\ol\l_q\!+\!2n\mu_q\!-\!b\mu_q)d_{p,q}^{2\mu,2n\mu+\l-\mu}
+(\ol\l_q\!+\!2n\mu_q\!-\!\mu_q\!+\!2b\mu_q)d_{p,q}^{-\mu,2n\mu+\l}
\vs{4pt}\hfill\cr
-(\ol\l_q\!+\!2n\mu_q\!+\!2b\mu_q)d_{p,q}^{-\mu,2(n+1)\mu+\l}
\!-\!(\ol\l_q\!+\!2n\mu_q\!+\!2\mu_q\!-\!b\mu_q)d_{p,q}^{2\mu,2n\mu+\l}
\!+\!3\mu_q d_{p,q}^{\mu,2n\mu+\l}\!=\!0.
\hfill\cr}
\eqno(3.42)$$
Taking the limit of $n\rar\infty$, using (3.41) (we assume that conditions
(3.29), (3.39)
also hold for $2\mu,-\mu$),
by noting that the sum of the second term and the third term vanishes
and so does the last term, we obtain
$$
{(\ol\l_{p''}-\mu_{p''})(\ol\l_{p''}+\mu_{p''})\over
(\ol\l_q-\mu_q+2b\mu_q)}d^{2\mu,\l-\mu}_{p,q}
-{\ol\l_{p''}(\ol\l_{p''}+2\mu_{p''})\over
(\ol\l_q+2b\mu_q)}d^{2\mu,\l}_{p,q}=0.
\eqno(3.43)$$
In (3.20), setting $r,\mu,\nu$ to be $q,2\mu,-\mu$ respectively, using
(3.31), (3.38), (3.43), we deduce
$$
d^{2\mu,\l}_{p,q}=
{(\ol\l_{p''}+2\mu_{p''})(\ol\l_q+b\mu_q)\over
2(\ol\l_{p''}+\mu_{p''})(\ol\l_q+2b\mu_q)}d^{\mu,\l}_{p,q}.
\eqno(3.44)$$
Using (3.44) in (3.43), we obtain
$$
{(\ol\l_{p''}-\mu_{p''})(\ol\l_{p''}+\mu_{p''})^2(\ol\l_q+(b-1)\mu_q)\over
2\ol\l_{p''}(\ol\l_q+(2b-1)\mu_q)^2}d^{\mu,\l-\mu}_{p,q}
-{\ol\l_{p''}(\ol\l_{p''}+2\mu_{p''})^2(\ol\l_q+b\mu_q)\over
2(\ol\l_{p''}+\mu_{p''})(\ol\l_q+2b\mu_q)^2}d^{\mu,\l}_{p,q}=0
\eqno(3.45)$$
Observe that the determinant of the coefficients of
$d^{\mu,\l-\mu}_{p,q},d^{\mu,\l}_{p,q}$ in (3.34), (3.45)
is a nonzero rational function on $\mu$
(when $\l$ is fixed), denote it by $\D(\mu)$. Thus we can choose
$\mu\in\G$ such that $\D(\mu)\ne0$, then we obtain that
$d_{p,q}^{\mu,\l}=0$.
Note that any element $\mu\in\G$ can be expressed as a sum of two elements
$\mu=\mu'+\mu''$ such that $\mu',\mu'',\mu'-\mu''$ all satisfy
conditions (3.29), (3.39) and $\si(\mu'),\si(\mu''),\D(\mu'),\D(\mu'')\ne0$.
Then by taking $r=q,
\mu=\mu',\nu=\mu''$ in (3.20), we obtain that $d_{p,q}^{\mu,\l}=0$ for all
$\mu\in\G$. If $\l$ does not satisfy (3.28), then using (3.20) again, we
can again deduce that $d_{p,q}^{\mu,\l}=0$.
Then (3.30) and (3.21) show that $e_{p,q}^{\mu,\l}=0$.
This proves Claim 2.
\par
{\bf Claim 3}. For $p\in\ol{1,\ell_1+\ell_2},\,q\in\ol{1,\ell},\,
\mu,\l\in\G,\,\vec i\in\Z_+^{\ell_1+\ell_2}$, we have
$$
\matrix{
(t^{1_{[p]}}\ptl_q)v_{\l,\vec i}\hfill&\!\!\!\!=&\!\!\!\!
(\a_q+\l_q)v_{\l,\vec i+1_{[p]}}
+(i_q+\d_{p,q}b)v_{\l,\vec i+1_{[p]}-1_{[q]}},
\vs{4pt}\hfill\cr
(x^\mu\ptl_q)v_{\l,\vec i+1_{[p]}}\hfill&\!\!\!\!=&\!\!\!\!
(\a_q+\l_q+b\mu_q)v_{\mu+\l,\vec i+1_{[p]}}
+(i_q+\d_{p,q})v_{\mu+\l,\vec i+1_{[p]}-1_{[q]}},
\vs{4pt}\hfill\cr
(x^\mu t^{1_{[p]}}\ptl_q)v_{\l,\vec i}\hfill&\!\!\!\!=&\!\!\!\!
(\a_q+\l_q+b\mu_q)v_{\mu+\l,\vec i+1_{[p]}}
+(i_q+\d_{p,q}b)v_{\mu+\l,\vec i+1_{[p]}-1_{[q]}}.
\hfill\cr}
\eqno\matrix{(3.46)\vs{4pt}\cr(3.47)\vs{4pt}\cr(3.48)\cr}$$
\par
Suppose $\vec i\ne0$. We shall use induction on $\vec i$.
Let $p'=p'_{\vec i},\,p''=p''_\l$.
Suppose $q=p''$. If $p\le p'$, then $p'_{\vec i+1_{[p]}}=p$ and so
(3.46) is the definition (3.9). If $p>p'$, then using (3.9), we obtain that
$(t^{1_{[p]}}\ptl_{p''})v_{\l,\vec i}$ is equal to
$$
\matrix{
(\a_{p''}+\l_{p''})^{-1}
((\d_{p',p''}t^{1_{[p]}}\ptl_{p''}-\d_{p,p''}t^{1_{[p']}}\ptl_{p''})
v_{\l,\vec i-1_{[p']}}
\vs{4pt}\hfill\cr
+t^{1_{[p']}}\ptl_{p''}
((\a_{p''}+\l_{p''})v_{\l,\vec i+1_{[p]}-1_{[p']}}
+(i_{p''}-\d_{p',p''}+\d_{p,p''}b)v_{\l,\vec i+1_{[p]}-1_{[p']}-1_{[p'']}})
\vs{4pt}\hfill\cr
-(i_{p''}-\d_{p',p''}+\d_{p',p''}b)
(t^{1_{[p]}}\ptl_{p''})v_{\l,\vec i-1_{[p'']}}),
\hfill\cr}
\eqno(3.49)$$
which is equal to the right-hand side of (3.46) by inductive assumption and
by noting that $\vec i+1_{[p]}-1_{[p']}<\vec i$ since $p>p'$ (cf. (3.1)).
If $q\ne p''$, we can prove (3.46) similarly.
Using induction, as in the proofs of (3.16), (3.17), we see that (3.47), (3.48)
hold up to some elements in $V_{\a+\mu+\l}^{(0)}$. Now the proof is exactly
analogous to that of Claim 2.
\par
{\bf Claim 4}. $B=\{v_{\mu,\vec i}\,|\,(\mu,\vec i{\ssc\,})\in
\G\times\Z_+^{\ell_1+\ell_2}\}$ is an $\F$-basis of $V$.
\par
By (3.10), it is straightforward to verify that $B$ is $\F$-independent.
So, it remains to prove that $V$ is spanned by $B$.
For any $\mu\in\G,\vec i\in\Z_+^{\ell_1+\ell_2}$, set
$$
\matrix{\dis
V_{\a+\mu}^{[\vec i]}=\{v\in V_{\a+\mu}\,|\,
\prod_{p=1}^\ell (\ptl_p-(\a_p+\mu_p))^{j_p}v
=0,\,\forall\,\vec j>\vec i\},\ \
V_{\a+\mu}^{(\vec i{\ssc\,})}=\bigcup_{\vec j<\vec i}V_{\a+\mu}^{[\vec j]},
\vs{4pt}\hfill\cr\dis
V^{[\vec i]}=\bigoplus_{\mu\in\G}V_{\a+\mu}^{[\vec i]},\ \
V^{(\vec i)}=\bigoplus_{\mu\in\G}V_{\a+\mu}^{(\vec i)}.
\hfill\cr}
\eqno(3.50)$$
Then
$$V_{\a+\mu}=\bigcup_{\vec j\in\sZ_+^{\ell_1+\ell_2}}V_{\a+\mu}^{[\vec j]}
\mbox{ and }
v_{\mu,\vec i}\in V_{\a+\mu}^{[\vec i]}\bs V_{\a+\mu}^{(\vec i)},
\ \ \mu\in\G,\ \vec i\in\Z_+^{\ell_1+\ell_2},
\eqno(3.51)$$
by (3.10) and induction on $\vec i$.
On the contrary, suppose that $V$ cannot be spanned by $B$ and that
 $\vec i$ is the minimal element in $\Z_+^{\ell_1+\ell_2}$
such that there exists $v\in V_{\a+\mu}^{[\vec i]}\bs V_{\a+\mu}^{(\vec i)}$
such that $v$ is not in the space spanned by $B$.
Then $\vec i\ne0$. Let $q$ be the minimal index
with $i_q\ne0$. Then
$v'=(\ptl_q-(\a_q+\mu_q))v\in V_{\a+_\mu}^{[\vec i-1_{[q]}]}$.
By the assumption, $v'$ is in the space spanned by $B$. Thus
$$
v''=v'-c{\ssc\,}v_{\mu,\vec i-1_{[q]}}\in
V_{\a+\mu}^{(\vec i-1_{[q]})},
\eqno(3.52)$$
for some $c\in\F$.
Let $w=v-c{\ssc\,}i_q^{-1}v_{\mu,\vec i}$.
Given any $\vec j\ge\vec i$. If $j_q\ne0$, then
$\vec j-1_{[q]}\ge\vec i-1_{[q]}$, and by (3.10), (3.51),
we have
$$
\prod_{p=1}^\ell(\ptl_p-(\a_p+\mu_p))^{j_p}w=
\prod_{p=1}^\ell(\ptl_p-(\a_p+\mu_p))^{j_p-\d_{p,q}}v''=0.
\eqno(3.53)$$
On the other hand, if $j_q=0$, then $\vec j>\vec i$, and so by (3.9), (3.50),
we have
$$
\prod_{p=1}^\ell(\ptl_p-(\a_p+\mu_p))^{j_p}w=
\prod_{p=1}^\ell(\ptl_p-(\a_p+\mu_p))^{j_p}v-
c{\ssc\,}i_q^{-1}\prod_{p=1}^\ell(\ptl_p-(\a_p+\mu_p))^{j_p}v_{\mu,\vec i}
=0.
\eqno(3.54)$$
(3.53), (3.54) show that $w\in V_{\a+\mu}^{(\vec i)}$
, and so $w$ is in the space spanned by $B$, a
contradiction with that $v$ is not in the space spanned by $B$. Thus Claim 4
follows.
\par
Now using Claim 3, we can obtain (2.12) by induction on $\vec i$. This
proves Theorem 2.3(1) in this case.
\par
{\it Case 2}. $\a\in\G$.
\par
Then $V(\a)=V(0)$. Thus we can suppose $\a=0$.
First assume that $V=V_0$. Then $(x^\mu t^{\vec i}\ptl_q)V\subset
V_\mu=0$ for all $\mu\in\G\bs\{0\}$, and since $W$ is generated by
$x^\mu t^{\vec i}\ptl_q$ with $\mu\ne0$, we have $WV=0$.
Since $V$ is indecomposable,
we obtain that $V=\F v_{0,0}$ must be the trivial submodule of $A_{0,0}$.
Thus assume that there exists $\mu\in\G\bs\{0\}$ such that $V_\mu\ne\{0\}$.
Then by [SZ1], $V_\mu\ne\{0\}$ for all $\mu\in\G\bs\{0\}$. Suppose
$V_0=\{0\}$. Then $V^{(0)}=\oplus_{\mu\in\G\bs\{0\}}V_{\mu}^{(0)}
\cong A'_{0,0}$ as $W(0,0,\ell_2+\ell_3,\G)$-modules,
and
$$
-\mu_p v_{\mu,0}=(x^{2\mu}\ptl_p)v_{-\mu,0}=
(x^\mu\ptl_p)(x^\mu t^{1_{[p]}}\ptl_p)v_{-\mu,0}
-(x^\mu t^{1_{[p]}})(x^\mu\ptl_p)v_{-\mu,0}=0,
\eqno(3.55)$$
for all $\mu\in\G,p\in\ol{1,\ell_1+\ell_2}$, which shows that $\ell_2=0$.
Then by assumption, $\ell_1,\ell_3\ge1$. We have
$$
-\mu_\ell v_{\mu,0}=(x^{2\mu}\ptl_\ell) v_{-\mu,0}=
([x^\mu\ptl_1,x^\mu t^{1_{[1]}}\ptl_\ell]
-[x^\mu t^{1_{[1]}}\ptl_1,x^\mu\ptl_\ell])v_{-\mu,0}=0,
\eqno(3.56)$$
for all $\mu\in\G,$ a contradiction since $\G$ is nondegenerate.
Thus $V_0\ne\{0\}$.
Now as in Case 1, we obtain that (2.12) holds for all $\mu,\nu\in\G$ with
$\nu,\mu+\nu\ne0$ (with suitable $b$).
Keeping this in mind will help our understanding in the following discussion.
If $V^{(0)}$ is a
$W(0,0,\ell_2+\ell_3,\G)$-module of the intermediate series of type $A_{0,b}$,
then exactly as in Case 1, we can obtain Theorem 2.3(1).
The only difference
is that we shall define $v_{0,\vec i}$ to satisfy
$$
(t^{2_{[p']}}\ptl_{p'})v_{0,\vec i-1_{[p']}}=
(i_{p'}+1)v_{0,\vec i},
\eqno(3.57)$$
(cf. (3.9)), where $p'=p'_{\vec i}$.
\par
It remains to consider that the $W(0,0,\ell_2+\ell_3,\G)$-module $V^{(0)}$
is one of the following four cases: $A(\b),B(\b),
A'_{0,0}\oplus\F v_{0,0},A'_{0,1}\oplus\F v_{0,0}$.
First suppose $V^{(0)}=A(\b)$.
Define $v_{0,\vec i}$ as in (3.57).
Let $p\in\ol{1,\ell_1+\ell_2},\,q\in\ol{\ell_1+1,\ell},\,\mu\in\G\bs\{0\}$.
By applying
${\rm ad\ssc\,}\ptl_r$ to $(x^\mu\ptl_q)v_{-\mu,1_{[p]}}$, using
$(x^\mu\ptl_q)v_{-\mu,0}=0$ by (2.10), we see that
$(x^\mu\ptl_q)v_{-\mu,1_{[p]}}\in V_0^{(0)}$, i.e.,
$(x^\mu\ptl_q)v_{-\mu,1_{[p]}}=c_{p,q,\mu}v_{0,0}$ for some
$c_{p,q,\mu}\in\F$.
Applying ${\rm ad\ssc\,}x^\nu\ptl_q$ with $\nu\ne0$ to it, using (2.10),
we obtain
$$
\matrix{
c_{p,q,\mu}(\nu_q+\b_q)v_{\nu,0}=
(\mu_q-\nu_q)(x^{\mu+\nu}\ptl_q)v_{-\mu,1_{[p]}}
+x^\mu\ptl_q(-\mu_qv_{\nu-\mu,1_{[p]}}+\d_{p,q}v_{\nu-\mu,0})
\vs{4pt}\hfill\cr
=(\mu_q-\nu_q)(-\mu_qv_{\nu,1_{[p]}}+\d_{p,q}v_{\nu,0})
-\mu_q((\nu_q-\mu_q)v_{\nu,1_{[p]}}+\d_{p,q}v_{\nu,0})
+\d_{p,q}(\nu_q-\mu_q)v_{\nu,0}
\vs{4pt}\hfill\cr
=-\d_{p,q}\mu_qv_{\nu,0},
\hfill\cr}
\eqno(3.58)$$
for all $\nu\in\G\bs\{0\}$. This shows that $c_{p,q,\mu}=0$ and so
$(x^\mu\ptl_q)v_{-\mu,1_{[p]}}=0$.
Applying ${\rm ad\ssc\,}x^\nu\ptl_r$ with $\nu,\mu+\nu\ne0,\,
r\in\ol{\ell_1+1,\ell}$
to this, we obtain
$$
\matrix{
0\!\!\!\!&=\!\!\!\!&
x^{\mu+\nu}(\mu_r\ptl_q-\nu_q\ptl_r)v_{-\mu,1_{[p]}}
+x^\mu\ptl_q(-\mu_rv_{\nu-\mu,1_{[p]}}+\d_{p,r}v_{\nu-\mu,0})
\vs{4pt}\hfill\cr&
=\!\!\!\!&
\mu_r(-\mu_qv_{\nu,1_{[p]}}+\d_{p,q}v_{\nu,0})
-\nu_q(-\mu_rv_{\nu,1_{[p]}}+\d_{p,r}v_{\nu,0})
\vs{4pt}\hfill\cr&&
-\mu_r(\nu_q-\mu_q)v_{\nu,1_{[p]}}-\d_{p,q}\mu_rv_{\nu,0}+
\d_{p,r}(\nu_q-\mu_q)v_{\nu,0})
\vs{4pt}\hfill\cr&
=\!\!\!\!&
-\d_{p,r}\mu_qv_{\nu,0},
\hfill\cr}
\eqno(3.59)$$
which is impossible. Thus the case $V^{(0)}=A(\b)$ does not occur.
\par
Next suppose $V^{(0)}=B(\b)$.
Let $p\in\ol{1,\ell_1+\ell_2},q\in\ol{\ell_1+1,\ell},p\ne q,\,
\mu,\nu,\mu\pm\nu\in\G\bs\{0\}$. We shall calculate
$[x^\nu\ptl_q,x^{-\nu}\ptl_p](x^\mu t^{1_{[p]}}\ptl_q)v_{-\mu,0}$. It is equal
to
$$
\matrix{\!\!\!
-(\nu_q\ptl_p\!+\!\nu_p\ptl_q)
(x^\mu t^{1_{[p]}}\ptl_q)v_{-\mu,0}
\!\!\!\!\!&=\!-([\nu_q\ptl_p\!+\!\nu_p\ptl_q,x^\mu t^{1_{[p]}}\ptl_q]
\!+\!(x^\mu t^{1_{[p]}}\ptl_q)(\nu_q\ptl_p\!+\!\nu_p\ptl_q))
v_{-\mu,0}
\vs{4pt}\hfill\cr&
=\nu_q(\mu_q+\b_q)v_{0,0}.
\hfill\cr}
\eqno(3.60)$$
On the other hand, it is equal to
$$
\matrix{
x^\nu\ptl_q([x^{-\nu}\ptl_p,x^\mu t^{1_{[p]}}\ptl_q]v_{-\mu,0}
-(\mu_p+\nu_p)(x^\mu t^{1_{[p]}}\ptl_q)v_{-\mu-\nu,0})
\vs{4pt}\hfill\cr
-x^{-\nu}\ptl_p([x^\nu\ptl_q,x^\mu t^{1_{[p]}}\ptl_q]v_{-\mu,0}
+(\nu_q-\mu_q)(x^\mu t^{1_{[p]}}\ptl_q)v_{\nu-\mu,0}).
\hfil\cr}
\eqno(3.61)$$
The first term of (3.61) is equal to
$$
\matrix{
x^\nu\ptl_q((x^{\mu-\nu}\ptl_q
+x^{\mu-\nu}t^{1_{[p]}}(\mu_p\ptl_q+\nu_q\ptl_p))v_{-\mu,0}
-(\mu_p+\nu_p)(-\nu_q)v_{-\nu,1_{[p]}})
\vs{4pt}\hfill\cr
=-\nu_q(x^\nu\ptl_q)v_{-\nu,0}
+x^\nu\ptl_q(-2\mu_p\nu_qv_{-\nu,1_{[p]}}+\nu_qv_{-\nu,0})
+\nu_q(\mu_p+\nu_p)(x^\nu\ptl_q)v_{-\nu,1_{[p]}}
\vs{4pt}\hfill\cr
=\nu_q(\nu_q+\b_q)v_{0,0}
-2\mu_p\nu_q(x^\nu\ptl_q)v_{-\nu,1_{[p]}}
-\nu_q(\nu_q+\b_q)v_{0,0}
+\nu_q(\mu_p+\nu_p)(x^\nu\ptl_q)v_{-\nu,1_{[p]}}
\vs{4pt}\hfill\cr
=\nu_q(\nu_p-\mu_p)(x^\nu\ptl_q)v_{-\nu,1_{[p]}}.
\hfill\cr}
\eqno(3.62)$$
Similarly, the second term of (3.61) is equal to
$$
-x^{-\nu}\ptl_p(
(\mu_q-\nu_q)(x^{\mu+\nu}t^{1_{[p]}}\ptl_q)v_{-\mu,0}
+(\nu_q-\mu_q)\nu_qv_{\nu,1_{[p]}})=0.
\eqno(3.63)$$
(3.60)-(3.63) show that
we obtain that
$(\mu_q+\b_q)v_{0,0}=
(\nu_p-\mu_p)(x^\nu\ptl_q)v_{-\nu,1_{[p]}}$ for all
$\mu,\nu,\nu\pm\nu\in\G\bs\{0\}$ with $\nu_q\ne0$. This is impossible.
Thus $V^{(0)}=B(\b)$ does not occur either.
\par
Now suppose $V^{(0)}=A'_{0,0}\oplus\F v_{0,0}$.
Let $\mu\in\G,q\in\ol{1,\ell},p\in\ol{1,\ell_1+\ell_2}$.
Then $(x^\mu\ptl_q)v_{-\mu,0}=0$.
Applying ${\rm ad\ssc\,}\ptl_r,\,r\in\ol{1,\ell}$
to $u=(x^\mu t^{1_{[p]}}\ptl_q)v_{-\mu,0}$, we see that
$u\in V_0^{(0)}$,
and so $(x^\mu\ptl_q)u=0$. Thus we have (3.55), which leads to
$\ell_2=0$. So $\ell_1,\ell_3\ge1$.
Let
$p\in\ol{1,\ell_1},\,q\in\ol{\ell_1+1,\ell},\,
\mu,\nu,\nu\pm\mu\in\G\bs\{0\}$.
We have $(x^\mu t^{1_{[p]}}\ptl_q) v_{-\mu,0}=c_{p,q,\mu}v_{0,0}$
for some $c_{p,q,\mu}\in\F$.
Let $\ptl=\sum_{r=1}^{\ell_1}a_r\ptl_r\in\DD_1$ (cf. (2.2)).
Then
$$
c_{p,q,\mu}(x^\nu\ptl)v_{0,0}\!=\!
[x^\nu\ptl,x^\mu t^{1_{[p]}}\ptl_q]v_{-\mu,0}
\!+\!(x^\mu t^{1_{[p]}}\ptl_q)(x^\nu\ptl)v_{-\mu,0}
\!=\!a_px^{\mu+\nu}\ptl_qv_{-\mu,0}\!=\!-a_p\mu_qv_{\nu,0},
\eqno(3.64)$$
for all $p\in\ol{1,\ell_1},\,q\in\ol{\ell_1+1,\ell},\,
\mu,\nu,\nu\pm\mu\in\G\bs\{0\}$. This shows that
$(x^\nu\ptl)v_{0,0}\ne0$ for all $\nu\ne0,\ptl\in\DD_1\bs\{0\}$.
Thus the linear map $\tau:\DD_1\rar V_\nu^{(0)}$ defined by
$\tau(\ptl)=(x^\nu\ptl)v_{0,0}$ is injective
But ${\rm dim}_{\sF}V_{\nu}^{(0)}=1$, so $\ell_1=1$.
Thus we can simply denote $\vec i$ just by $i\in\Z_+$.
(3.64) shows that
$(x^\nu\ptl_1)v_{0,0}$ does not depend on $\nu$. By rescaling $v_{0,0}$ if
necessary, we can suppose $(x^\nu\ptl_1)v_{0,0}=v_{\nu,0}$ for all
$\nu\in\G\bs\{0\}$. Now if we denote
$v'_{0,0}=0,v'_{0,i}=v_{0,i-1}$ and $v'_{\nu,i}=v_{\nu,i},\nu\ne0$, then
as in case 1, we can prove
$$
(x^\mu t^i\ptl_q)v'_{\nu,j}=
\left\{\matrix{
jv'_{\mu+\nu,i+j-1}\hfill&\mbox{if }q=1,
\vs{4pt}\hfill\cr
\nu_qv'_{\mu+\nu,i+j}\hfill&\mbox{if }q>1.
}\right.
\eqno(3.65)$$
This shows that $V=A'_{0,0}$ as $W(1,0,\ell-1,\G)$-modules.
\par
Finally suppose $V^{(0)}=A'_{0,1}\oplus\F v_{0,0}$. Then as above, we have
$\ell_2=0$, and we also have (3.64) except the last equality which should
now be $a_p\nu_qv_{\nu,0}$. Since we can always take $\nu$ such that
$\nu_q\ne\nu_{q'}$ for all $q,q'\in\ol{\ell_1,\ell_1},q\ne q'$. This
shows that $\ell_1=\ell_3=1$. Thus $A'_{0,1}=A'_{0,0}$ as
$W(0,0,1,\G)$-modules and this case becomes the previous case.
This completes the proof of Theorem 2.3(1).
\par
The proof of the rest of Theorem 2.3 is standard (see also [Z1,SL]).
\par\
\vs{-7pt}\par
\cl{\bf References}
\par
\ni\hi3ex\ha1
[BXZ] G.~Benkart, X.~Xu and K.~Zhao, ``Classical Lie superalgebras
over the simple associative algebras,'' preprint 2000.
\par
\ni\hi3ex\ha1
[DZ1] D.~Z.~Dokovic, K.~Zhao, ``Derivations, isomorphisms, and
second cohomology of generalized Witt algebras,''
{\it Trans.~of Amer.~Math.~Soc.}, {\bf350}\,(1998), 643-664.
\par\ni\hi3ex\ha1
[DZ2] D.~Z.~Dokovic, K.~Zhao, ``Generalized Cartan type $W$ Lie algebras
in characteristic zero,'' {\it J.~Alg.}, {\bf195}\,(1997), 170-210.
\par\ni\hi3ex\ha1
[E1] S.~Eswara Rao, ``Irreducible representations of the Lie algebra of
the diffeomorphisms of a $d$-dimensional torus,'' {\it J.~Alg.},
{\bf182}\,(1996), 401-421.
\par\ni\hi3ex\ha1
[E2] S.~Eswara Rao, ``Representations of Witt algebras,''
{\it Publ.~Res.~Inst.~Math.~Sci.}, {\bf30}\,(1994), 191-201.
\par\ni\hi3ex\ha1
[K] V.~G.~Kac, {\it Vertex Algebras for Beginners,}
University Lectures Series, Vol.~{\bf10}, AMS. Providence RI, 1996.
\par\ni\hi3ex\ha1
[Kaw] N.~Kawamoto, ``Generalizations of Witt algebras over a field of
characteristic zero,'' {\it Hiroshima Math.~J.}, {\bf 16}\,(1986), 417-462.
\par\ni\hi3ex\ha1
[M] O.~Mathieu, ``Classification of Harish-Chandra modules over the
 Virasoro Lie algebra,'' {\it Invent.~Math.}, {\bf 107}\,(1992), 225-234.
\par\ni\hi3ex\ha1
[P] D.~P.~Passman, ``Simple Lie algebras of Witt type,'' {\it J.~Alg.},
{\bf 206}\,(1998), 682-692.
\par\ni\hi3ex\ha1
[PZ] J.~Patera, H.~Zassenhaus, ``The higher rank Virasoro algebras,''
 {\it Comm.~Math.~Phys.}, {\bf 136}\,(1991), 1-14.
\par\ni\hi3ex\ha1
[S1] Y.~Su, ``A classification of indecomposable $sl_2(\C)$-modules
 and a conjecture of Kac on irreducible modules over the Virasoro
 algebra,'' {\it J. Alg.}, {\bf 161}\,(1993), 33-46.
\par\ni\hi3ex\ha1
[S2] Y.~Su, ``Harish-Chandra modules of the intermediate series over
 the high rank Virasoro algebras and high rank super-Virasoro
 algebras,'' {\it J.~Math.~Phys.}, {\bf 35}\,(1994), 2013-2023.
\par\ni\hi3ex\ha1
[S3] Y.~Su, ``Classification of Harish-Chandra modules over the
 super-Virasoro algebras,''$\ssc\!$ {\it Comm.~Alg.}, {\bf 23}\,(1995), 3653-3675.
\par\ni\hi3ex\ha1
[S4] Y.~Su, ``Simple modules over the high rank Virasoro algebras,''
 {\it Comm.~Alg.}, in press.
\par\ni\hi3ex\ha1
[SX] Y.~Su, X.~Xu, ``Structure of divergence-free Lie algebras,'' to appear.
\par\ni\hi3ex\ha1
[SXZ] Y.~Su, X.~Xu and H.~Zhang, ``Derivation-simple algebras and the
structures of Lie algebras of Witt type,'' {\it J.~Alg.}, in press.
\par\ni\hi3ex\ha1
[SZ1] Y.~Su, K.~Zhao, ``Generalized Virasoro and super-Virasoro algebras
and modules of the intermediate series,''
to appear.
\par\ni\hi3ex\ha1
[SZ2] Y.~Su, K.~Zhao, ``Second cohomology group of generalized Witt type
Lie algebras and certain representations,''
to appear.
\par\ni\hi3ex\ha1
[SZ3] Y.~Su, K.~Zhao, ``Simple algebras of Weyl type,''
to appear in {\it Sciences in China}.
\par\ni\hi3ex\ha1
[SZ4] Y.~Su, K.~Zhao, ``Isomorphism classes and automorphism groups
of algebras of Weyl type,'' to appear in {\it Sciences in China}.
\par\ni\hi3ex\ha1
[SZh] Y.~Su, L.~Zhu ``Non-graded Virasoro and super-Virasoro
algebras and modules of the intermediate series,'' to appear
\par\ni\hi3ex\ha1
[X1] X.~Xu, ``New generalized simple Lie algebras of Cartan type over a field
with characteristic 0,'' {\it J.~Alg.}, {\bf224}\,(2000), 23-58.
\par\ni\hi3ex\ha1
[X2] X.~Xu, {\it Algebraic Theory of Hamiltonian Superoperators,}
monograph, to appear.
\par\ni\hi3ex\ha1
[Z1] K.~Zhao, ``Harish-Chandra modules over generalized Witt algebras,''
to appear.
\par\ni\hi3ex\ha1
[Z2] K.~Zhao, ``Isomorphisms between generalized Cartan type $W$ Lie
algebras in characteristic zero,'' {\it Canadian J.~Math.}, {\bf 50}\,(1998),
210-224.
\end{document}